\documentclass[11pt]{amsart}
   \usepackage{amsthm,amsmath,amssymb,amscd,graphicx,enumerate, stmaryrd,xspace,verbatim, epic, eepic,color,url}

\usepackage{eurosym}

\usepackage[all]{xypic}
\SelectTips{cm}{}
\usepackage{pdfsync}
 
   \newtheorem{thm}{Theorem}[subsection]
      \newtheorem*{thm*}{Theorem}

   \newtheorem{prop}[thm] {Proposition}

   \newtheorem*{conjecture*}{Conjecture}

{\theoremstyle{definition}
 
          \newtheorem*{exercise*}{Exercise}

  \newtheorem{remark} [thm]{Remark}

\newcommand{\N}{{\mathbb{N}}}

\newcommand{\RR}{{\mathbb{R}}}

\newcommand{\Z}{{\mathbb{Z}}}

\newcommand{\cC}{{\mathcal C}}

\renewcommand{\cH}{{\mathcal H}}

\renewcommand{\cL}{{\mathcal L}}
\newcommand{\cM}{{\mathcal M}}

\newcommand{\cP}{{\mathcal P}}

\newcommand{\cX}{{\mathcal X}}

 \newcommand{\md}{{\underline{d}}}

  \newcommand{\mt}{{\underline{t}}}

\def\<{\langle}
\def\>{\rangle}

\newcommand{\oO}{{\overline{O}}}

\newcommand{\Spec}{\operatorname{Spec}}

 \newcommand{\Sing}{{\operatorname{Sing}}}

\newcommand{\Pic}{{\operatorname{{Pic}}}}

\newcommand{\Aut}{{\operatorname{Aut}}}

\newcommand{\trop}{{\operatorname{trop}}}

\newcommand{\an}{{\operatorname{an}}}

\newcommand{\double}{\genfrac..{0pt}1
{\raise -2pt\hbox{$\scriptstyle\longrightarrow$}}{\raise 4pt\hbox
{$\scriptstyle\longrightarrow$}}}

 \newcommand{\la}{\longrightarrow}
\newcommand{\ha}{\hookrightarrow}

\newcommand{\ov}{\overline}

 \newcommand{\Mgbst}{\ov{\mathcal{M}_g}}
  \newcommand{\Mgban}{\Mgb^{\operatorname{an}}}
    \newcommand{\Mgnban}{\Mgnb^{\operatorname{an}}}
 
\def\Mgb{\overline{M}_g}
 
 \newcommand{\Mgnb}{\ov{M}_{g,n}}
  \newcommand{\Mgnbst}{{\ov{\cM}}_{g,n}}
  
  \newcommand{\Mgnst}{{\mathcal{M}_{g,n}}}
\def\sPdgb{\overline{\mathcal{P}}_{g}^d}
\def\sPdg{{\mathcal{P}_{g}^d}}

\newcommand{\PX}{\overline{P}^{g-1}_X}
\newcommand{\PXd}{\overline{P}^{d}_X}

\newcommand{\Pdgb}{\overline{P}^{d}_g}

\newcommand{\sPgg}{\overline{\cP}^{g}_g}
\newcommand{\PXg}{\overline{P}^{g}_X}

\def\agale{\overline{A}^{{\text{mod}}}_g}

\newcommand{\Mgtnb}{\ov{M}_{g,n}^{\rm trop}}

   \newcommand{\oOG} {\overline{\mathcal{OP}}}

\newcommand{\oOP}{{\overline{\mathcal{OP}}}}

 \newcommand{\Sg}{{\mathcal {SG}}_g} 
  \newcommand{\Sgn}{{\mathcal {SG}}_{g,n}}

\newcommand{\HM}{\mathcal{A}_{\bullet}}

\begin{document}

\bibliographystyle{plain}

\title{Recursive combinatorial aspects of  compactified moduli spaces}

\author[]{Lucia Caporaso}

  \address[Caporaso]{Dipartimento di Matematica e Fisica\\
 Universit\`{a} Roma Tre\\Largo San Leonardo Murialdo \\I-00146 Roma\\  Italy }\email{caporaso@mat.uniroma3.it}

 \date{\today}




\maketitle

\tableofcontents

 \section{Introduction}
  
In recent   years an interesting connection has been established
between some  moduli spaces of algebro-geometric objects 
(e.g. algebraic stable curves)
and some moduli spaces of polyhedral
 objects (e.g.   tropical curves).

In loose words, this connection expresses the   Berkovich skeleton  of a given algebro-geometric moduli space  as the moduli
  space of the skeleta of the objects parametrized by  the given    space;
it has been proved to hold   in two important   cases: the moduli space  of stable curves 
and the moduli space  of admissible covers.  Partial results are known in  other cases.

 This connection  relies on the study of the
 boundary of the algebro-geometric moduli spaces and on its
  recursive, combinatorial   properties, some of which have  been long   known  and    are now
   viewed from a  new perspective.

We will describe some  results in this area by focusing on    the  moduli spaces of curves, line bundles on curves (i.e. Jacobians),  and coverings of curves.
As we said, it was clear for a long time that combinatorial aspects play a significant  role in compactifying moduli spaces. A notable
example    is the structure of the N\'eron model of the Jacobian of a curve. Of course, N\'eron models are not compact, but they 
 are a first step towards compactifying  Jacobians. We begin the paper  by describing N\'eron models, their combinatorial
properties, and the recursive structure of their  compactification.
Then we turn to moduli spaces of curves and illustrate the connection introduced at the beginning.
We return to compactified Jacobians in the last section  and describe some recent partial results.

 \section{Compactified Jacobians and N\'eron models}
\label{NJsec}
 \subsection{Jacobians, Picard schemes and N\'eron models}

 \label{Neronsec}
  Let $X$ be a connected, reduced,  projective curve of genus $g$ over an algebraically closed field $k$. 
The set of isomorphism classes of line bundles of degree $0$ on every irreducible component of $X$ is an algebraic group of dimension $g$,
 denoted by $J_X$, and called the Jacobian  of $X$.
 More exactly, $J_X$ is the moduli space for line bundles of multidegree $(0,0,\ldots, 0)$ on $X$.
 
$J_X$ is projective, hence an abelian variety, if $X$ is smooth   but not if $X$ is singular  
(with   exceptions).  
Constructing  compactifications for $J_X$ is a classical problem which     can be stated as follows.
Let
$$
X\ha \cX\stackrel{f}{\la} \Spec R
$$
be a family of curves over the spectrum of a  discrete valuation ring $R$,  with $X$ as special fiber. Assume the generic fiber, $\cX_K$, to be a nonsingular curve over the quotient field, $K$, of $R$, and  let
$J_{\cX_K}\to \Spec K$ be its Jacobian.  The   problem is to find ``good" models of $J_{\cX_K}$ over $\Spec R$, where ``good" means:    (a) projective,     (b)  with a  moduli interpretation,  (c) such that the special fiber depends only on $X$  and   not  
on   the family $f$.
We shall refer to such models, and to their special fiber, as 
{\it compactified Jacobians}.

We shall begin by looking for models that satisfy (b) and (c) together with the weaker requirement of being separated,  rather than projective.

A natural  model  of $J_{\cX_K}$  is the relative Jacobian (a group scheme over $R$)
\begin{equation}
\label{jac}
 J_{\cX/R} \la \Spec R.
\end{equation}
This is    a   separated model for $J_{\cX_K}$, but it  does not have   satisfactory moduli properties.
Indeed, suppose we have a line bundle $\cL$ over $\cX$ having relative degree $0$ on the fibers of $f$. Then
we have a {\it moduli morphism}  $\mu_{\cL_K}:\Spec K\to J_{\cX_K}$ whose image corresponds to the restriction of $\cL$ to $\cX_K$.
For a   model   to have good moduli properties we want the map $\mu_{\cL_K}$ to extend 
to a morphism from $\Spec R$ to the model (so that the image of  the special point is determined by the restriction of $\cL$ to the special fiber). 
For the relative Jacobian     this requirement easily fails. 

Another   natural model is the relative degree-$0$ Picard scheme
\begin{equation}
\label{pic}
 \Pic^0_{\cX/R}\la \Spec R.
 \end{equation}
Now,  this has  a good moduli interpretation,
as any moduli  map $\mu_{\cL_K}$ will certainly extend to a map $\mu_R:\Spec R \to  \Pic^0_{\cX/R}$,
and $\mu_R$ is itself a moduli map. The problem now  is that the extension $\mu_R$ may fail to be unique, that is,  
\eqref{pic} is not a separated morphism, in general.

We have thus two natural  models, one   separated with bad moduli properties, the other with good moduli properties but not separated. Does there exist a compromise with better behaviour than both of them?
An answer to this  question comes from the theory of N\'eron models.

Our moduli requirement above (i.e. the existence of a unique extension for   maps of the form $\mu_{\cL_K}$)
is a special case of the  mapping property satisfied by N\'eron models,
whose existence has been established in \cite{Neron}. 
We state the following special case of  N\'eron's  famous theorem, using the above notation.
\begin{thm}[N\'eron] 
\label{N} 
 There exists a smooth and separated group scheme of finite type 
 $$
N(J_{\cX_K})\la \Spec R,
$$
whose generic fiber  is $J_{\cX_K}$,   satisfying the following   mapping property:

Let  $Y_R\to \Spec R$ be   smooth, $Y_K$   its generic fiber, and     $\mu_K: Y_K\to J_{\cX_K}$   a morphism.
Then $\mu_K$ extends uniquely to a morphism $\mu_R:Y_R\to N(J_{\cX_K})$.
\end{thm}
\begin{remark}
The N\'eron model does not commute with ramified base change, therefore  it is
 a   separated model for $J_{\cX_R}$ with good mapping properties, but not   functorial ones.
\end{remark}

The     relation  of the N\'eron model  with
the relative Jacobian and the relative Picard scheme has ben established by Raynaud who proved,  in   \cite{Raynaud}, 
that 
 $N(J_{\cX_K})\to \Spec R$ is the maximal separated quotient of $ \Pic^0_{\cX/R}\to \Spec R$.

 \subsection{Combinatorics of N\'eron models}
 \label{Ncomb}
We shall   turn to the  geometric  structure of our N\'eron models.
We assume from now on that the curve $X$ is nodal, write $X=\cup_{v\in V}C_v$ for its irreducible components, and consider the dual graph, $G_X$, of $X$:
\begin{equation}
 \label{GX}
G_X:=(V, E=\Sing(X), V\stackrel{h}{\la} \Z)
\end{equation}
with $h(v)=g(C_v)$, where $g(C_v)$ is the   genus of  the normalization, $C_v^{\nu}$.
The genus of $G_X$ is   the same as the arithmetic genus of $X$, i.e.
$$
g(X)=g(G_X)=b_1(G_X)+\sum_{v\in V}h(v).
$$

It is well known that  $J_X$ fits  into an exact sequence of algebraic groups:
\begin{equation}
 \label{es}
0\la (k^*) ^{b_1(G_X)} \la J_X \la \prod_{v\in V}J_{C_v^{\nu}}\la 0.
\end{equation}
 
Denote by $N_X$ the special fiber of $
N(J_{\cX_K})\to \Spec R
$. 
We shall   state some results of Raynaud, \cite{Raynaud}, and Oda-Seshadri, \cite{OS},   which establish that $N_X$ is a disjoint union of copies of the Jacobian of $X$ 
indexed by a  combinatorial invariant of the curve.
 
To do that we need some combinatorial preliminaries. Fix an  orientation on $G=G_X$ (whose choice is irrelevant),
let 
 $C_0(G, \Z)$ and $C_1(G, \Z)$ be the standard groups of $i$-chains, generated over $\Z$ by $V$ if $i=0$,
and by $E$  if $i=1$. Next, let 
 $\partial:C_1(G, \Z)\to C_0(G, \Z)$ be  the usual boundary 
 (mapping an edge $e$ oriented from $u$ to $v$ to $u-v$), and 
 $\delta:C_0(G, \Z)\to C_1(G, \Z)$ the coboundary (mapping a vertex $v$ to $\sum e_v^+-\sum e_v^-$
 where the first sum is over all edges originating from $v$ and the second over  all edges ending at $v$).
 
 Now we can state the following.
 
\begin{prop}[Raynaud, Oda-Seshadri]
\label{ROS}
 Let $\cX\to \Spec R$ have nodal special fiber and regular total space.
 Then the special fiber, $N_X$, of the N\'eron model $N(J_{\cX_K})\la \Spec R$
is a union of copies of $J_X$ as follows
\begin{equation}
 \label{NX}
N_X\cong \sqcup_{i\in \Phi_X}(J_X)_i
\end{equation} where $\Phi_X$  is the following finite group 
$$
\Phi_G\cong \frac{\partial \delta C_0(G, \Z)}{\partial C_1(G \Z)}. 
$$
In particular, the number of irreducible components of $N_X$ equals the number of spanning trees of $G_X$.
\end{prop}
 The last  claim follows from   Kirchhoff-Trent, or Kirchhoff matrix,Theorem.
 
 \subsection{Compactifying N\'eron models}
 \label{compNer}
 From now on we   apply the notation introduced in Proposition~\ref{ROS} and for any connected nodal curve $X$ we   denote by $N_X$ the special fiber of the N\'eron model of the Jacobian
 associated to a family  $\cX\to \Spec R$ with $\cX$ regular.

 If $X$ is a singular curve with $b_1(G_X)\neq 0$ (i.e. $X$ not of ``compact type"), then 
  the N\'eron model $N(J_{\cX_K})\to \Spec R$
 is not projective, as its special fiber is not projective  by  the exact sequence \eqref{es}.
  
Now,   the   Picard scheme has good  moduli properties 
 and   the N\'eron  model is its maximal separated quotient.
We   introduce  a terminology to distinguish compactified Jacobians 
which also compactify the N\'eron model.
 
 A   compactified Jacobian  $ 
\ov{P}\la \Spec R
 $
 is said to be 
 of {\it N\'eron type} (or a {\it N\'eron compactified Jacobian}) if  its special fiber, written $\ov{P}_X$,  contains $N_X$ as a dense open subset.
 We shall also say that $\ov{P}_X$ is of N\'eron type.
  
  As we shall see,  N\'eron compactified Jacobians   do exist, but there exist also interesting compactified Jacobians
  not of N\'eron type.
  
The notion  originates from \cite{OS} and \cite{Cneron}, although  the terminology was introduced later, in \cite{Ctype}.
Oda and Seshadri, in \cite{OS}, treated the case of a fixed singular  curve $X$ (rather than a  family of curves),
and constructed compactified Jacobians in this less general sense.
 They nonetheless  established the link with N\'eron models and  constructed compactified Jacobians
whose smooth locus is isomorphic to  $N_X$.
This was extended to families of curves later. In \cite{Cthesis}  and \cite{Cneron} a class of compactified Jacobians  of N\'eron type was proved to exist and to form a family over the moduli space of stable curves, $\Mgbst$.
Such families, denoted by  $
   \psi_d:\sPdgb \la \Mgbst 
$, are indexed by the integers $d$ such that 
\begin{equation}
 \label{MR}
 (d-g+1,2g-2)=1.
\end{equation}
For any stable curve $X$ the fiber of $\psi_d$ over $X$, written 
$ 
\PXd
$, 
contains $N_X$ as a dense open subset equal to its smooth locus, and it is thus of N\'eron type.

 In a similar vein, in \cite{MV} and in \cite{MRV}  other N\'eron  compactified Jacobians were found among the ones constructed by Esteves in \cite{Esteves},
 and called ``fine" compactified Jacobians. 
The word ``fine" is quite  appropriate, as   all known N\'eron  compactified Jacobians
 are  as fine a moduli space as they can be, i.e.  they   admit a  universal 
 (or ``Poincar\'e") line bundle.


We   now 
show how  N\'eron Jacobians are  
  recursive compactifications of N\'eron models.
 First, for a  connected graph $G$, we 
introduce the following  $$
  \cC(G):=\{S\subset E(G):   \  G-S  \text{ is connected}\},
  $$  
with partial  order given by reverse inclusion. 
The       maximal  element of $ \cC(G)$  is $\emptyset$, and the   
  minimal elements   are the    $S\subset E$ such that $G-S$ is a spanning tree. 
Moreover, $\cC(G)$
 is a graded poset with respect to the     rank function
 $S\mapsto g(G-S).$ 

 If $G$ is the dual graph of the curve $X$ then $S\in \cC(G)$ is a set of nodes of $X$. We denote by
  $X_S^{\nu}$ the desingularization of $X$ at $S$, so that $X_S^{\nu}$ is a connected nodal curve of genus $g(G-S)$.
  Recall that $N _{X_S^{\nu}}$ denotes the special fiber of the N\'eron model of its Jacobian.
 The following   follows from  \cite[Thm. 7.9]{Cneron}.
 
\begin{thm} 
\label{Nstratthm}
 Let $\PXd$ be a N\'eron compactified Jacobian. Then   
\begin{equation}
 \label{Nstrat}
\PXd=\bigsqcup_{S\in \cC(G)} N_S
\end{equation}
 with $N_S
\cong N _{X_S^{\nu}}$ for every $S\in \cC(G)$. Moreover \eqref{Nstrat} is a \emph{graded stratification}, i.e.  the following hold. 
\begin{enumerate}
\item
 $
N_S\cap \overline{N_{S'}}\neq \emptyset \Leftrightarrow  N_S\subset \overline{N_{S'}} 
 \Leftrightarrow S'\geq S.
 $
 \item
 $N_S$ is locally closed of pure dimension $g(G-S)$.
 \item
 The following is a rank-function on $\cC(G)$
 $$
\cC(G) \la \Z;\quad \quad S\mapsto \dim N_S.  
 $$
\end{enumerate}
\end{thm}
\begin{remark}
 The set of strata of minimal dimension in \eqref{Nstrat} are N\'eron models of curves 
 whose dual graph is a spanning tree of $G_X$, hence they are irreducible.
 By Proposition~\ref{ROS}, the number of such   strata is equal to the number of irreducible components of $\PXd$.
 \end{remark}
If $X$ is not of compact type     the strata of \eqref{Nstrat} are not all connected. Hence
one naturally asks whether the stratification can be refined so as to have connected strata. We will answer this question later in the paper.
 
\begin{remark}
\label{rec}
 The theorem exhibits   the compactification of the N\'eron model of $X$
in terms of the N\'eron models of partial normalizations of $X$.
This phenomenon is an instance of what seems to be a widespread  recursive  behaviour for compactified moduli spaces.
Namely, to compactify a space (e.g. $N_X$) one  adds at the boundary 
the analogous spaces   associated to simpler objects (e.g. $N_{X'}$ with $X'$ a connected partial normalization of $X$).
Other examples of this recursive pattern will appear later in the paper.
 \end{remark}

The concept of  ``graded stratification" used in the Theorem will appear again, so we now define it in general.

A {\it{graded stratification}}  of an algebraic variety,  or a stack,  $M$   by a poset  $\cP$ is 
  a  partition 
$ 
 M=\bigsqcup_{p\in \cP}M_p
$ 
such that the following hold for every $p,p'\in \cP$.
\begin{enumerate}
 \item
$ M_p\cap \overline{M_{p'}}\neq \emptyset \Leftrightarrow  M_p\subset \overline{M_{p'}} 
 \Leftrightarrow p'\geq p.
 $
 \item
$M_p$ is    equidimensional and locally closed.
 
\item
The map  from $\cP$ to $\N$ sending $p$ to $\dim M_p$
is a rank function on $\cP$.
\end{enumerate}

\section{Moduli   of curves: tropicalization and analytification}
\label{Mgsec}
\subsection{Moduli spaces of algebraic and tropical curves}
Let $\Mgnst$ be the moduli space of smooth curves of genus $g$ with $n$ marked points.
Assume $2g-2+n>0$  so that  $\Mgnst$ is never empty.
It is well known that $\Mgnst$ is not projective (unless $g=0, n=3$)  and   is compactified
by  the moduli space of Deligne-Mumford stable curves, $\Mgnbst$; see  \cite{DM}, \cite{Knudsen2}, \cite{Knudsen3}, \cite{Gieseker}.

We will describe in \eqref{stratMg} a   stratification of $\Mgnbst$ which is recursive in the sense of   Remark~\ref{rec},
that is,  the boundary strata are expressed 
in terms  of  simpler moduli spaces $\cM_{g',n'}$. 
In this case ``simpler"means   $g\leq g'$ and of smaller dimension,
 i.e. $n'<n+3(g-g')$. These boundary strata will be described in \eqref{MG}.

Let $\Sgn$ be the  poset of  stable graphs of genus $g$ with $n$ legs, with the following partial order:
$$
G_2\geq G_1 \quad  \text{ if  } \quad   G_2=G_1/S
\quad\text{  for some } \quad S\subset E(G_1) 
$$
 where $G_1/S$ is the graph obtained by contracting every edge of $S$ to a vertex.
 By the very definition  (introduced in \cite{BMV}), edge-contraction preserves $g$ and $n$,  and it is easily seen to preserve
stability.  

To a stable curve, $X$, of genus $g$ with $n$ legs there corresponds a dual graph
$G_X\in \Sgn$. With respect to what we defined   in \eqref{GX}  the  only  new piece of data are the legs of $G_X$  which correspond  to the $n$ marked points. Recall that the weight, $h(v)$,  of a vertex $v$ is the geometric genus of the corresponding component of $X$.
We denote by $\deg(v)$ the degree (or valency) of  $v$. 

For $G\in \Sgn$  we denote by $\cM_G$ the moduli stack of   curves whose dual graph is isomorphic to $G$.
We have (see \cite[Prop. 3.4.1]{ACP})
\begin{equation}
 \label{MG} 
\cM_G=\Bigr[\bigr(\prod_{v\in V(G)}\cM_{h(v), \deg(v)}\bigl)/\Aut(G)\Bigl].
\end{equation}
With this notation and   the   terminology at the end of Section~\ref{NJsec}  we state:
\begin{prop}
\label{propalg} 
The following is a graded stratification of $\Mgnbst$ 
\begin{equation}
 \label{stratMg} 
 \Mgnbst= \bigsqcup_{G\in \Sg}\cM_G.
\end{equation}
 \end{prop}

One goal of the above   descriptive result (whose proof, in \cite{Chbk},  is not hard thanks to our consolidated knowledge of $\Mgbst$)
 is to highlight  the similarities 
of $\Mgnbst$ with the moduli space of extended (abstract) tropical curves, $\Mgtnb$, as we are going to show.
First of all,  an  extended tropical curve is a metric graph, i.e.
a graph $G$   whose edges are assigned a length,
$\ell:E\to \RR_{>0}\cup \{\infty\}$. We denote a tropical curve as follows
$$
\Gamma= (G; \ell)=(V,E,w; \ell).
$$
The genus of $\Gamma$ is the genus of $G$.
The word ``extended"   refers to the fact that we include   edges of infinite length.
In fact, abstract   tropical curves  are originally defined (in \cite{MZ} and \cite{BMV}) as compact   spaces, so their  edges   have finite length.
We allow edges of infinite lengths     to obtain a 
 compact moduli space which  
 is  the ``tropicalization" of the moduli space of stable curves.

An extended tropical curve with $n$ marked points is a metric graph as above with the addition of a set of legs of the underlying graph.

There is a natural equivalence relation on tropical curves such that in every equivalence class there is a unique (up to isomorphism)  curve whose underlying graph is stable.

The moduli space $\Mgtnb$ parametrizes extended tropical curves of genus $g$ with $n$ marked points up to this equivalence relation. 
Its first construction is due to Mikhalkin, for $g=0$, and to Brannetti-Melo-Viviani for compact tropical curves in any  genus.
The following statement summarises results from \cite{MIK5}, \cite{BMV} and \cite{Chbk}.

\begin{thm}
\label{thmtrop}
 The moduli space of extended tropical curves, $\Mgtnb$, is a compact and normal topological space of dimension $3g-3+n$.
 It admits a graded stratification   
\begin{equation}
 \label{stratMgt} 
 \Mgtnb= \bigsqcup_{G\in \Sgn^*}\overline{M}_G^{\trop}
\end{equation}
 where $\overline{M}_G^{\trop}$ is the locus of  curves having $G$ as underlying graph and $\Sgn^*$ is the  poset  dual to $\Sgn$ (i.e. with reverse partial order).
 \end{thm}

 \subsection{Skeleta and tropicalizations}
\label{sktrop}
Proposition~\ref{propalg} and Theorem~\ref{thmtrop}
show that $ \Mgtnb$ has a graded stratification dual to that of $\Mgnbst$.
Therefore we   ask whether there exists
some deeper relation between $\Mgnbst$ and $\Mgtnb$.

A positive answer can be given 
through the theory of analytifications of algebraic schemes developed in   \cite{Berkovich}, and through its connections to tropical geometry;
see \cite{MS} and \cite{Payne}.
Let us introduce the   space
$\Mgnban$, the   analytification of $\Mgnbst$ in the sense of Berkovich. Recall that a point in $\Mgnban$ corresponds, up to base change,
 to a stable curve over an algebraically closed  field $K$   complete with respect to a non-archimedean valuation;
as $\Mgnb$ is projective, this is the same as a stable curve over  the ring of integers of $K$.

From the general theory (see also \cite{Thuillier} and \cite{ACP})
we have that for every space with a toroidal structure, like the stack $\Mgnbst$ (with toroidal structure given by  its boundary), one associates the {\it Berkovich skeleton} which is a generalised, extended cone complex  onto which the analytification
retracts. We write $\overline{\Sigma}(\Mgnbst)$ for the Berkovich skeleton of $\Mgnbst$
and $ 
 \Mgnban \stackrel{\rho}{\la} \overline{\Sigma}(\Mgnbst) 
 $ 
for the retraction.

The connection between tropical geometry and Berkovich theory originates from the fact that
Berkovich skeleta can be viewed as tropicalizations of algebraic varieties. 
As we are going to see, the picture for curves is quite clear, whereas other interesting  situations 
(some treated later in the present paper) are
still open to investigations;
we refer to  \cite{GRW}, for recent progress for higher dimensional varieties. 

Now,
tropical curves can be viewed as tropicalizations (or skeleta) of curves over algebraically closed, complete, non-archimedean fields.
 This is clarified by the following statement, combining   results of  \cite{Viviani}, \cite{Tyomkin} \cite{BPR1}, \cite{BPR2}
(for the first part)
and of  \cite{ACP} (for the second part).

\begin{thm}
\label{ACP}
There exists a    {\emph {tropicalization}}  map $$ \trop: \Mgnban\to  \Mgtnb$$
that  sends the class of a stable curve  over the (algebraically closed, complete, non-archimedean)  field  $K$  to its     {\emph {skeleton}}.

There is a
 natural isomorphism,  $\overline{\Sigma}(\Mgnbst)\cong \Mgtnb$,
through which the above map  factors as follows
 $$
 \trop: \Mgnban \stackrel{\rho}{\la}  \overline{\Sigma}(\Mgnbst) \stackrel{\cong}{\la} \Mgtnb.
 $$
 \end{thm}
We need to define the tropicalization map   and explain  the word ``skeleton".
As we said, a point in $\Mgnban$ is a class of stable curves,
$\cX\to \Spec R$, where $R$ is   the   valuation ring of $K$; let   $X$ be its special fiber. Then the image of this point via the map $\trop$ is the tropical curve $(G; \ell)$, 
where $G$ is the dual graph of $X$ and, for every node $e\in E(G)$, the value $\ell(e)$ is determined by the local geometry of $\cX$ at $e$ as measured by the given valuation.  
Such a tropical curve $(G;\ell)$ is called the {\it skeleton}, of the stable curve $\cX\to \Spec R$
 or   of the stable curve $\cX_K$  over   $K$.


Concluding in loose words: the skeleton of $\Mgnbst$ is the moduli space of skeleta of    stable curves.

The structure of the proof of the above theorem is  such that it may apply to other situations. In fact 
it has  been   applied by Cavalieri-Markwig-Ranganathan to another remarkable case, the compactification of the Hurwitz spaces, as we shall now explain.

\subsection{Algebraic and tropical admissible covers}
\label{sktropAC}
Consider 
the moduli space of admissible covers,  $\ov{\cH}_{\bullet}$, which here (as in \cite{CMR}), is the normal compactification
of
 the classical   space of Hurwitz covers ${\cH}_{\bullet}$.
We use the subscript 
 $``\bullet"$ to simplify the notation 
 needed to express the usual discrete invariants. 
Indeed, a more precise notation 
 would be  ${\cH}_{\bullet}=\cH_{g,h}\bigr(\underline{\pi}\bigl)$ for the Hurwitz space
 parametrizing  degree-$d$ covers of  a smooth curve of genus $h$ by a smooth curve of genus $g$ 
with exactly $b$ branch points with ramification profile prescribed by a set of $b$ partitions of $d$,
written $\underline{\pi}=(\pi_1,\ldots, \pi_b)$.

 In \cite{CMR} the authors define {\it tropical admissible covers},
 construct their moduli space $\ov{H}^{\trop}_{\bullet}$ and establish an analogue to Theorem~\ref{ACP}. 
 We shall now outline
 the  procedure and give some details.

The first step is to associate  a dual combinatorial entity to the algebro-geometric one.
Indeed, one associates to an  admissible cover
 a map of   graphs, which we   call the  {\it dual graph cover}.
The  set of all dual graph covers is endowed with a poset structure   by means of  edge-contractions (similarly   to 
  the poset set of stable graphs).
 We denote by 
  $\HM$  this poset, as its objects  can be viewed as  {\it admissible} maps of graphs. 
  For any $\Theta\in \HM$ we denote by  $\cH_{\Theta}$ the locus in $\ov{\cH}_{\bullet}$
  of   admissible covers whose dual graph map is $\Theta$.

The second step is to enrich the dual combinatorial entity with  a tropical, or polyhedral, structure. 
In \cite{CMR}     tropical admissible covers  are defined by metrizing, in a suitable way, dual maps of graphs.
The  moduli space of tropical admissible covers is denoted by $\ov{H}^{\trop}_{\bullet}$, the bar over $H$ indicates that they are ``extended", i.e. edge-lengths are allowed to be infinite.
 For   $\Theta\in \HM$ the   stratum parametrizing tropical admissible  covers having $\Theta$ as underlying graph map    
is denoted by $\ov{H}^{\trop}_{\Theta}$ and shown to be the quotient of an extended  real  cone.

The third and last step is to use analytification  and tropicalization  to establish an explicit  link between the algebraic and the tropical
moduli space.
Indeed, essentially  by construction, we have  dual stratifications
$$
\ov{\cH} _{\bullet}=\bigsqcup_{\Theta\in \HM}\cH_{\Theta}\quad\quad  \text{ and }\quad  \quad 
\ov{H}^{\trop}_{\bullet}=\bigsqcup_{\Theta\in \HM^*}\ov{H}^{\trop}_{\Theta}.
$$
This duality can   be read   from the following Theorem.
  \begin{thm}[Cavalieri-Markwig-Ranganathan]
  \label{CMR}
 There is a tropicalization map
 $ \trop:\ov{\cH}_{\bullet}^{\an}\to  \ov{H}^{\trop}_{\bullet} $
which  factors as follows
 $$
 \trop:\ov{\cH} _{\bullet}^{\an}\stackrel{\rho}{\la}  \ov{\Sigma}(\ov{\cH} _{\bullet}) \la  \ov{H}^{\trop}_{\bullet}.
 $$
\end{thm}
The map $\ov{\cH} _{\bullet}^{\an}\stackrel{\rho}{\la}  \overline{\Sigma}(\ov{\cH} _{\bullet})$ is the retraction of $\ov{\cH} _{\bullet}^{\an}$ onto its Berkovich skeleton, as described in subsection~\ref{sktrop}.
This result is compatible with the analogous one  for $\Mgbst$   through the canonical forgetful maps
from $\ov{\cH}_{\bullet}$ to the moduli spaces of stable curves; see \cite[Thm. 4]{CMR}.

  What about other moduli spaces?
Are Theorems~\ref{ACP} and \ref{CMR} part of some  general picture where skeleta of algebraic moduli spaces
(e.g.  the skeleton of $\Mgnbst$)
can be described as moduli spaces for combinatorial entities
(e.g. tropical curves) which are skeleta of   the   objects (e.g. stable curves)  parametrized by the algebraic moduli spaces?
 As we saw, the first step is to identify  a suitable partially ordered set of combinatorial objects to associate
to the algebro-geometric ones.

In the next section we shall look at the theory of compactified Jacobians from this point of view.

\section{Compactified  Jacobians}
\subsection{Compactifying  Jacobians over $\Mgb$}
Let us go back to  compactify  Jacobians of curves  and, with the discussion of the previous section in mind,   
  approach the problem from the point of view of the  moduli theory of stable curves.
Consider the universal Jacobian over the moduli space of  smooth curves   and  look for  a compactification of it
over $\Mgbst$ satisfying the requirements we discussed earlier. 

For reasons that will be clear later, it is convenient to extend our considerations  to Jacobians of all degree.
For a curve  $X=\cup_{v\in V}C_v$ and any multidegree $\md\in \Z^V$
we write $ \Pic^{\md}(X)$ for the moduli space of 
 line bundles of multidegree $\md$. Now, $J_X$ is identified with $\Pic^{(0,\ldots,0)}(X)$ and  we have non-canonical isomorphisms  $J_X\cong  \Pic^{\md}(X)$.

The  universal degree-$d$ Jacobian over $\cM_g$  is a   morphism 
$
 \sPdg \to \cM_g,
$
whose fiber over the point parametrizing a smooth curve $X$ is $\Pic^d(X)$.
We want to  construct a 
compactification of $\sPdg$ over $\Mgbst$ by  a projective morphism
  $ 
\sPdgb \to \Mgbst
  $ 
such that $\sPdgb$
has a moduli description. We shall refer to such a space as a     compactified universal degree-$d$ Jacobian.

From  this perspective  there is a natural  approach to the problem,    namely imitate 
the construction of $\Mgbst$ itself. Recall that 
the moduli scheme, $\Mgb$, 
of the stack $\Mgbst$
was   constructed    by   Gieseker in \cite{Gieseker}
as the GIT-quotient of the Hilbert scheme of $n$-canonically embedded curves ($n\gg  0$). 
The  moduli stack $\Mgbst$ is   the quotient stack associated to this quotient.


Now, as the Hilbert scheme of $n$-canonically embedded   curves of genus $g$  
  has such a beautiful GIT-quotient,  why shouldn't the Hilbert scheme
  of all  projective curves of  degree $d\gg 0$ and genus $g$  have a beautiful GIT-quotient?
  And why shouldn't this quotient be a candidate for a compactification of the
  universal degree-$d$ Jacobian?
  Indeed, this is   what happens, 
  and the GIT quotient of this Hilbert scheme is our compactified   universal degree-$d$ Jacobian, $\Pdgb$.
The corresponding quotient stack is denoted by $\sPdgb$.

Now,   as $d$ varies, the spaces $\sPdgb$ are not isomorphic to one another,
indeed their fibers over certain singular curves are not even birational to one another.
Again, we see  the phenomenon  (appearing already  in \cite{OS})
that  non isomorphic compactifications of the Jacobian of a fixed singular curve exist.
 In the present case   the various models depend on the degree $d$.

As for the  basic properties of the spaces $\sPdgb$ as $d$ varies, there is a special
set of degrees, namely those  such that
  \eqref{MR} holds,
such that  $\Pdgb$  is a geometric GIT-quotient and has good moduli properties,
so  that its points correspond  to geometric  objects up to a certain  equivalence relation.
Moreover,  the natural (projective) 
 morphism     
 $ 
  \psi_d:\sPdgb \la \Mgbst
 $ 
  is a strongly representable map of Deligne-Mumford stacks.  
  As we said in subsection~\ref{compNer} in this case all fibers of $\psi_d$ are N\'eron compactified Jacobians.

Now, \eqref{MR} fails
if $d=g-1$  whereas it holds 
 if $d=g$.
 We shall concentrate on these two cases, interesting for different reasons,
and    give a combinatorial analysis of the compactified Jacobian.

\subsection{Compactified Jacobians in degree $g-1$}
We begin by reviewing an idea of Beauville.
Let us fix   $g\geq 2$.
Among   degree-$d$ Jacobians, the case $d=g-1$ has been object of special interest
for its strong connections with the Theta divisor, the Schottky problem,
the Prym varieties; in particular, it has    been studied in \cite{Beauville}.

Let us approach  the problem of compactifying  the degree $(g-1)$-Jacobian of a curve $X$ with $G=(V,E,h)$ as dual graph.
We expect   a good  compactification to have   finitely many irreducible components
  and  each component  to parametrize  (at least away from the boundary)
line bundles on $X$ of a fixed multidegree $\md$ such that $|\md|=g-1$.
Now the question is to determine these ``special" multidegrees.
Consider the identity
$$
g-1=\sum_{v\in V}\bigr(h(v)-1\bigl)+|E|. 
$$
We can interpret the first term (i.e. the summation) as carrying the topological invariants of $X$, and the second term, $|E|$,
as carrying the combinatorial ones. Now, while the first term exhibits the single contribution of each vertex/component,
the second does not. So we may ask  how  to distribute the second term among the various vertices in a combinatorially meaningful way. A natural solution   is to consider an orientation, $O$, on $G$, and denote by 
$\mt^O_v$   the number of edges having $v$ as target. Then 
$ 
\sum_{v\in V}\mt^O_v=|E|.
$ 
Therefore, if we define  a multidegree $\md^O$  as follows 
$$
\md^O_v:=h(v)-1+\mt^O_v,
$$
for every $v\in V$, we have $|\md^O|=g-1$.
As  there are only finitely many orientations on a graph,   by the above rule we have picked a finite set 
of  special  multidegrees of degree $g-1$. More precisely, it may happen 
that two orientations, $O$ and $O'$,  give the same multidedegree; in such a case we say that $O$ and $O'$ are {\it equivalent}.
We denote by 
$ 
\oO(G)
$ 
the set of such equivalence classes of orientations on $G$.
Now, a closer look  reveals that $\oO(G)$ is still too big  for it to parametrize the irreducible components of a compactified Jacobian.
Indeed, from the discussion in Subsection~\ref{Ncomb} we expect the number of such components to be at most equal to
$|\Phi_G|$, whereas we have $|\oO(G)|>|\Phi_G|$ in general
(for example for a 2-cycle, as in the picture below).

So, to compactify the degree-$(g-1)$ Jacobian we must    distinguish a special type of orientations. These are called
{\it totally cyclic} orientations, defined as orientations such that any two vertices
in the same connected component of  $G$
lie in a directed cycle.
It tuns out that if two orientations  are equivalent,   one is totally cyclic  if the other one is.
The set of equivalence classes of totally cyclic orientations on $G$ is denoted by 
$ 
\oO^0(G).
$

In the picture below we have the four orientations on a 2-cycle. The first two are totally cyclic and equivalent. The last two are not totally cyclic.
    
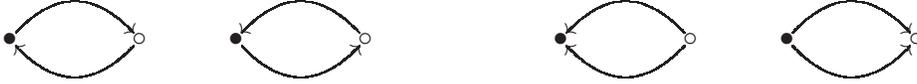
\begin{figure}[ht]
\begin{equation*}
$$\xymatrix@=.5pc{
 *{\bullet} \ar @{->} @/^1.2pc/[rrrr] \ar @{<-} @/_1.2pc/[rrrr]&&&& *{\circ} 
 &&&*{\bullet} \ar @{<-} @/^1.2pc/[rrrr] \ar @{->} @/_1.2pc/[rrrr] ]&&&& *{\circ} 
 &&&&&&*{\bullet} \ar @{<-} @/^1.2pc/[rrrr] \ar @{<-} @/_1.2pc/[rrrr] ]&&&& *{\circ}
 &&&*{\bullet} \ar @{->} @/^1.2pc/[rrrr] \ar @{->} @/_1.2pc/[rrrr] ]&&&& *{\circ}
 \\
}$$
\end{equation*}
\end{figure}

We shall adopt the convention that the empty orientation on the graph consisting of only one vertex and no edges
is totally cyclic.
We notice the following  facts.
\begin{remark}
\label{rk0}
\begin{enumerate}[(a)]
\item 
 $ \oO^0(G)$ is empty if and only if $G$ contains some bridge.  
 
\item
\label{rk0c} Assume $G$ connected. Then 
$| \oO^0(G)|\leq |\Phi_G|$ with equality if and only if $|V|=1$.
\end{enumerate}
\end{remark}

The set of all   orientations on the spanning subgraphs of a  graph $G$ admits a partial order as follows.
Let $O_{S_1}$ and $O_{S_2}$ be two orientations  on $G-S_1$ and $G-S_2$ respectively, where $S_i\subset E(G)$ for $i=1,2$.
We set
$O_{S_1}\leq O_{S_2}
$ 
if $G-S_1\subset G-S_2$ and if the restriction of $O_{S_2}$  to $G-S_1$ equals $O_{S_1}$.

This definition is compatible with the equivalence relation defined above, and hence the set of all equivalence classes of totally cyclic orientations on $G$ is a poset, which we shall denote as follows
$$
\oOP^0(G):= \bigsqcup_{S\subset E(G)} \oO^0(G-S).
$$

Finally, we are ready to exhibit  a graded stratification  of $\PX$ governed by totally cyclic orientations, by 
 rephrasing some  results in \cite{CV1}.
 \begin{prop}  
\label{PX}
Let $X$ be a stable curve of genus $g$ and $G$ its dual graph. Then the following is a graded stratification 
  \begin{equation}
 \label{strato}
\PX =\bigsqcup_{{\overline{O_S}}\in  \oOG^0(G)}P_X^{O_S},
\end{equation}
and we have     natural isomorphisms for every  ${\overline{O_S}}\in  \oOG^0(G)$
\begin{equation}
 \label{isoo}P_X^{O_S}\cong \Pic^{\md^{O_S}}(X^{\nu}_S).
 \end{equation}
 \end{prop}

The isomorphisms \eqref{isoo} exhibit  the  recursive behaviour 
described in Remark~\ref{rec}.
Indeed, $\Pic^{\md^{O_S}}(X^{\nu}_S)\cong  J_{X^{\nu}_S}$ and $\md^{O_S}$ is the multidegree associated to a totally cyclic orientation on $G-S$.
Hence   the boundary of the compactified degree-$(g-1)$ Jacobian of $X$
is stratified by  Jacobians  of  degree $(g(X')-1)$ of   partial normalizations, $X'$, of $X$.

\begin{remark}
If $X$ is reducible and not of compact type then 
$\PX$ is not of N\'eron type, by   Remark~\ref{rk0}\eqref{rk0c}. 
\end{remark}

A tropical counterpart of the stratification \eqref{strato} is  not   known to us.
 
The compactified Jacobians in degree $g-1$ have been proved especially useful in   various situations and in connection with the   Theta divisor, whose definition extends to these compactified Jacobians.
Among the applications, we recall that 
the pair given by this degree-$(g-1)$ Jacobian and its Theta divisor, $(\PX, \ov{\Theta}_X)$, is endowed with a natural group action of   $J_X$ and, as such,
forms a so-called {\it principally polarized stable semi-abelic  pair}. These pairs appear as boundary points in the 
compactification, $\agale$, of the moduli space of principally polarized abelian varieties constructed in \cite{Alexeev}.
Moreover, by  \cite{AlexeevTorelli}, they form the image of the compactified Torelli map 
and  we have
\begin{prop}
The extended Torelli  morphism
$ 
\ov{\tau}: \Mgb\to   \agale 
$ 
mapping a curve $X$ to $(\PX, \ov{\Theta}_X)$ is the moduli map associated to the family 
$\psi_{g-1}:\ov{P}^{g-1}_g\to\Mgb$.
\end{prop}
The combinatorial structure of $\PX$ described in Proposition~\ref{PX}   is heavily used in \cite{CV1} to study  the  fibers of $\ov{\tau}$.

 \subsection{Jacobians in degree $g$: N\'eron compactified Jacobians}
 We now adapt the considerations at the beginning of the previous sections to the case $d=g$.
 We have
 $ 
g=\sum_{v\in V}\bigr(h(v)-1\bigl)+|E| +1
$ 
 and, as before,  we want to  express the term $|E|+1$ in a combinatorially meaningful way.
Modifying what we did earlier, we now consider
orientations having one  bioriented edge. More precisely, 
a {\emph {1-orientation}} on a graph $G$ is the datum of a {\it {bioriented}} edge, $e$,
 and of 
  an orientation on $G-e$.   For any 1-orientation we have, of course, 
 $ 
|\md^O| =g(G). 
$ 
Just as in the previous subsection  we restricted to totally cyclic orientations, now we need to restrict to certain special 1-orientations,
namely   ``rooted"  orientations.
A 1-orientation with bioriented edge $e$ is said to be {\it rooted} (or $e${\it -rooted})
  if for every  vertex $v$ there exists a directed path from $e$ to $v$.
 
  As before,   two 1-orientations $O$ and $O'$ such that  $\md^O=\md^{O'}$  are defined to be equivalent,
  moreover one is rooted if the other one is.
   We   denote by $\oO^1(G)$ the set of equivalence classes of rooted $1$-orientations on $G$.
  
We decree the empty orientation on the graph consisting of only   vertex and no edges
to be  rooted.
Similarly to Remark~\ref{rk0} we have:
\begin{remark}
\label{rk1}
\begin{enumerate}[(a)]
\item 
 $ \oO^1(G)$ is not empty if and only if $G$ is connected.  
 \item
\label{rk1c}
$| \oO^1(G)|=|\Phi_G|$.
\end{enumerate}
\end{remark}

In the picture below we have all $e$-rooted orientations on a 4-cycle  with fixed bioriented edge $e$.
They correspond to the four elements in  $\oO^1(G)$.
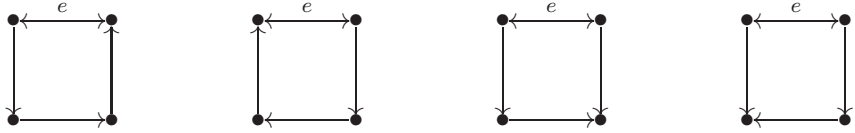
\begin{figure}[!htp]
$$\xymatrix@=1pc{
 &*{\bullet} \ar@{>}[dd]
\ar@{<->}[rr]^e& & *{\bullet}\ar@{<-}[dd]
&&&
*{\bullet} \ar@{<-}[dd]
\ar@{<->}[rr]^e& & *{\bullet}\ar@{>}[dd]
&&&
*{\bullet} \ar@{>}[dd]
\ar@{<->}[rr]^e& & *{\bullet}\ar@{>}[dd]
&&&*{\bullet} \ar@{>}[dd]
\ar@{<->}[rr]^e& & *{\bullet}\ar@{>}[dd] \\
&&&  &   &&&& &&&&& & &&&&&&&&&&&& &&&&&&&&&&  \\
&*{\bullet} \ar@{>}[rr]& &*{\bullet} &&&
*{\bullet} \ar@{<-}[rr]& &*{\bullet}&&& 
*{\bullet} \ar@{>}[rr]& &*{\bullet}&  &&
*{\bullet} \ar@{<-}[rr]& &*{\bullet}&&&& 
}$$
\label{four}
\end{figure}

In Theorem~\ref{Nstratthm} we exhibited a stratification indexed by the poset of connected spanning subgraphs,  $\cC(G)$.
The strata of that stratification are not connected, we shall now 
exhibit a finer stratification with connected strata.

Recall that rooted $1$-orientations exist only on connected graphs.
The poset of equivalence classes of rooted $1$-orientations on all spanning subgraphs of $G$
is written as follows \begin{equation}
 \label{OP1}
\oOP^1(G):= \bigsqcup_{S\in \cC(G)} \oO^1(G-S),
\end{equation}
with the   partial order   defined in the previous section.

The following result of  \cite{Christ} states that  $\PXg$  admits a recursive graded stratification     governed by rooted orientations.
 \begin{thm} [Christ]
\label{PXg}
Let $X$ be a stable curve of genus $g$ and $G$ its dual graph. Then $\PXg$   admits the following   graded stratification 
  \begin{equation}
 \label{strato1}
\PXg =\bigsqcup_{{\overline{O_S}}\in  \oOG^1(G)}P_X^{O_S},
\end{equation}
and we have     natural isomorphisms for every  ${\overline{O_S}}\in  \oOG^1(G)$
\begin{equation}
 \label{iso1}P_X^{O_S}\cong \Pic^{\md^{O_S}}(X^{\nu}_S).
 \end{equation}
 \end{thm}
 
\begin{remark}
 Comparing with the stratification of Theorem~\ref{Nstratthm} we have, using \eqref{OP1},
 that  \eqref{strato1} is a refinement of \eqref{Nstrat}, with connected strata.
\end{remark}

 In this case we do have a tropical counterpart. 
First, consistently with the dual   stratifications \eqref{stratMg} and \eqref{stratMgt}, 
a tropical counterpart of $X$ is a tropical curve $\Gamma$ whose underlying graph is the dual graph of $X$.

We now assume $\Gamma=(G,\ell)$ is compact (i.e. not extended).
The tropical curve $\Gamma$ has a Picard group $\Pic(\Gamma)=\sqcup \Pic^d(\Gamma)$,
and each connected component, $\Pic^d(\Gamma)$, is isomorphic to the same  $b_1(G)$-dimensional real torus; see \cite{MZ}.
In \cite{ABKS}, An-Baker-Kuperberg-Shokrieh    show  that  $\Pic^g(\Gamma)$  has an interesting  polyhedral decomposition indexed by ``break divisors" on $G$.
 The connection between break divisors and rooted 1-orientations is established, as a consequence of the results in \cite{ABKS}, in \cite{Christ},
 where
the following result is obtained.
  \begin{thm} 
\label{PXgtrop} Let $\Gamma=(G, \ell)$ be a compact tropical curve of genus $g$.
 Then $\Pic^g(\Gamma)$   admits the following   graded stratification 
  \begin{equation}
 \label{stratot}
\Pic^g(\Gamma)=\bigsqcup_{{\overline{O_S}}\in  \oOG^1(G)^*}\Sigma_{\Gamma}^{O_S}.
\end{equation}
 \end{thm}

 The stratification \eqref{stratot} is a non-trivial rephrasing of the polyhedral decomposition established in \cite{ABKS}. Such a  rephrasing  
is needed to establish the connection with the stratification
\eqref{strato1}.
 From \cite{ABKS} it follows that
the strata $\Sigma_{\Gamma}^{O_S}$   are the interiors of the faces of a polyhedral decomposition for $\Pic^g(\Gamma)$.

Presently, we do not know whether the duality between the stratifications  \eqref{stratot} and \eqref{strato1}
can be given an interpretation in terms of tropicalization and analytification,
similarly to the cases described in subsections~\ref{sktrop} and \ref{sktropAC}.


This problem is related to a result of \cite{BR}, which we will state in our  notation.
Using  Theorem~\ref{ACP}, let $\cX_K$ be a smooth curve,   $J_{\cX_K}$ its Jacobian,  and let $\Gamma=\trop([\cX_K])$,
where $[\cX_K]$ is the point of $\Mgban$ corresponding to $\cX_K$. Hence $\Gamma$ is a  compact tropical curve of genus $g$  
(compactness follows from    $\cX_K$  being smooth). With this   set-up,    \cite[Thm. 1.3]{BR} yields

\begin{thm}[Baker-Rabinoff]  
 $\Pic^g(\Gamma)\cong \Sigma(J_{\cX_K}^{\an})$.
\end{thm}

With  this   result  in mind,  a natural approach to  the problem   mentioned above  would be to study the relation
between $ \Sigma(J_{\cX_K}^{\an})$ and $ \PXg$. 
  
Finally, consider the    universal compactified Jacobian.  Results from \cite{Christ} indicate that an analogue of    Theorem~\ref{PXg}  should hold
 uniformly over $\Mgbst$, so that the universal compactified Jacobian $\sPgg$ can be given a graded stratification  compatible with that of $\Mgbst$.
We expect the same to hold for the universal compactified Jacobian  in degree $g-1$, with Proposition~\ref{PX} as starting point.
 
\bibliography{ICM}

\end{document}